\documentclass[10pt]{amsart}
\usepackage[utf8]{inputenc}
\usepackage[english,russian]{babel}

\usepackage{amsmath}
\usepackage{amssymb}
\usepackage{amsfonts}

\newtheorem{lemma}{Lemma}
\newtheorem{theorem}{Theorem}

\begin{document}
\renewcommand{\refname}{References}
\textbf{УДК 517.55+517.33}
\begin{center}
\textbf {ON SOME SHARP ESTIMATES OF TOEPLITZ OPERATOR IN
	SOME SPACES OF HARDY-LIZORKIN TYPE OF ANALYTIC
	FUNCTIONS IN THE POLYDISK 
 }\\

\vspace{8pt}

\textbf{R.F. Shamoyan}\\
\end{center}

\vspace{8pt}

\noindent{\small {\bf Abstract.} We provide some new sharp assertions on the action of Toeplitz
	$T_\varphi$ operator in new $F^{
	p,q}_\alpha$ type spaces of analytic functions of several complex
	variables extending previously known assertions proved by various authors.

\vspace{10pt}


\section{INTRODUCTION}

Let $U^n$ be the unit polydisk in $\mathbb{C}
^n, U^n = \{z \in \mathbb{C}^n : |z| < 1, j = 1, . . . ,n\}$. Let
further $H(U^n)$ be the space of all analytic functions in $U^n$. Let further also

$$ F^{p,q}_\alpha 
(U^n)=\left\{f \in H(U^n) : 
||f||^p_{F^{p,q}_\alpha}=\int\limits_{T^n}
\left(\int\limits_{I_n}
|f(r\xi)|^q
 (1-r)^{\alpha q-1}dr
  \right)^{p/q} d
  \xi < \infty
\right\},
$$

where $0 < p, q < \infty, \aleph > 0, T^n = \{|z_j | = 1, j = 1, . . ., n\}, I^n = (0, 1] \times · · · \times(0, 1],$
$dr = dr_1 ... dr_n,
d\xi = d\xi_1 ... d\xi_n,
 (1-r)^\alpha =\prod_{k=1}^n
 (1-r_k)^\alpha, r_k \in (0, 1)$ be the
holomorphic Lizorkin-Triebel space, (see, for example, [2], [3], [6], [9], [10]).

Note
for particular case $p = q$
 we have Bergman classical class, for $q = 2$ 
 we have so called Hardy-Lizorkin
 space $H^p_\beta$
for some $\beta$ that is,
 $H^p_\beta = 
 \left\{f\in H(U^n):D^\beta f\in H^p\right\},
 0 < p \leq \infty, \beta > 0,$ where $D^\beta$
is a fractional derivative of analytic $f$ function in $U^n$. Note (see definitions bellow) for this particular cases the action of $T_\phi$ classical
Toeplitz operator is well-studied in unit disk, unit ball and unit polydisk. We
study $T_\phi$ operators in more general $F^{p,q}_\alpha$ type spaces in the polydisk. 

Our main sharp result provide some criteria for symbol of $T_\phi$ to obtain boundednes of $T_\phi$ in
mentioned type analytic spaces.
We define classical Hardy space $H^p(U^n), 0 < p \leq \infty $ as follows (see also, for
example, [1] and [2]). Let

$$
H^p(U^n)=\left\{f \in H(U^n) : ||f||_{H^p}=\sup\limits_{r \in I^n}M_p(f,r)<\infty \right\}, 
$$
 and
 
 $$
 M_p(f,r)=\left(\int\limits_{T^n}|f(r\xi)|^p dm_n(\xi)\right)^{\frac{1}{p}},
 $$

where $r\xi=(r_1\xi_1, ... , r_n\xi_n )$ and $dm_n$ is a normalized Lebesgues measure on $T^n$,
$r_j \in (0, 1), j = 1, ... , n.$ Note $M_p(f, r) $function is growing function by each $r_j$
and we, for $p = \infty$ , obtain classical and well-studied class $H^\infty(U^n)$ of all bounded
analytic functions in $U^n$ (see for example [8] for this class of functions). Various
sharp results on action of Toeplitz operators can bee seen in papers of various
authors in various functional spaces in the unit ball and polydisk. We mention, for
example, the following papers [4], [5], where such type sharp results can be seen in
particular cases of $F^{p,q}_\alpha$ spaces namely in Bergman and in Hardy type spaces in the
unit ball and in the unit polydisk. We note similar type results in particular values
of parameters is well-known, (see, for example, [1], [2], [4], [5], [9]).
Such type sharp result on boundedness of Toeplitz operators also have various
applications (see, for example [1], [2], [5], [9]).
We remind the reader the standard definition of Toeplitz $T_h$ operators in the
unit polydisk.

Let $h \in L^1(T^n)$. Then we define Toeplitz $T_h$ operator as one integral operator

$$
(T_hf)(z)=\frac{1}{(2\pi)^n}\int\limits_{T_n}\frac{f(\xi_1,...,\xi_n)h(\xi_1,...,\xi_n)}{\prod\limits^{n}_{k=1}(1-\bar{\xi}_kz_k)}d\xi_1...d\xi_n,
$$

$k=1,...,n, z_k \in U.$

Note that we can easily show 
$F^{p,q}_\alpha$
general mixed norm analytic function spaces
in the unit polydisk are Banach spaces for all values of $p$ and $q$, if $min(p, q) > 1$
and they are complete metric spaces for all other values of $p$ and $q$.

We stress that behavior of the operators in the unit polydisk is substantially
different from the action of $T\phi$ operators in the unit ball in $\mathbb{C}^n$ (see [1], [4], [5], [9]
for example). Our intention to set criteria for the action of Toeplitz $T\phi$ operators
from $F^{p,q}_{\alpha,k} (U^n)$ into Bergman-Sobolev and Hardy-Lizorkin type spaces in the unit
polydisk, under the assumption that $\phi$is holomorphic,
 $\phi \in H(D^n)$ (some restriction
on symbol of Toeplitz operator).

We define some new function spaces in the polydisk for formulation of our main
result in the polydisk. Let further $dm_{2n}$ be the normalized Lebesgues measure in
$U^n, D^s, 0 < s \leq \infty$ be the fractional derivative of holomorphic $f$ function

$$
(D^sf)(z)=\sum_{|k|\geq 0}\frac{\Gamma(s+k+1)\Gamma(s+1)}{\Gamma(k+1)}a_kz^k,
$$

$a_kz^k=a_{k_1...k_m}z_1^{k_1}...z_m^{k_m}, f(z)=\sum_{|k\geq 0|}a_k^{z_k}, z\in U^n, \Gamma(\alpha+1)=\prod_{j=1}^m\Gamma(\alpha_j+1),$
$\alpha_j>-1, j=1,...,m.$

Note if $f \in H(U^n)$ then for any $s\in \mathbb{N}, D^sf\in H(U^n).$ 
Let also
 
$F^{p,q}_{\alpha,k} (U^n)=\left\{f\in H(U^n):||D^kf||_{F^{p,q}_\alpha}<\infty\right\}, 0<p,q,\alpha<\infty, k \in \mathbb{N}.$

$A^{s}_{\alpha,m} (U^n)=F^{s,s}_{\alpha/s,m}=$

$=\left\{f\in H(U^n):||f||^s_{A_{\alpha,m}}
\int\limits_{U^n}|(D^mf)(z)|^s(1-|z|)^{\alpha-1}
dm_{2n}(z) <\infty \right\},$

 $m\in \mathbb{N}, 0 < s, \alpha <\infty $
  (Bergman-Sobolev space). Let further

$$H^{s}_{m} (U^n)=\left\{f\in H(U^n):||D^mf||_{H_{s}}
 <\infty \right\}, m\in \mathbb{N}, 0 < s <\infty$$

be analytic Hardy-Lizorkin space in the unit polydisk $U^n.$

Note it can easily shown that these both scales of analytic function spaces in the
unit polydisk are Banach spaces for all values of $s, s \geq1$ and they are complete
metric spaces for other values of $s, s > 0$.
Throughout the paper, we write $C$ or $c$ (with or without lower indexes) to denote
a positive constant which might be different at each occurrence (even in a chain of
inequalities), but is independent of the functions or variables being discussed.

\section{Main results}

We prove our main sharp theorems in this section.

\begin{theorem}
	Let $0<\max (p,q)\leq s, 1<s<\infty, m,k \in \mathbb{N}, k=\alpha+\frac{1}{p}-\frac{1}{s}+m\left(1-\frac{1}{s}\right).$ Then $T_\varphi$ operator is bounded operator from $F^{p,q}_{\alpha,k} (U^n)$ into $A^{s}_{m,m} (U^n)$ if and only if $\varphi \in H^\infty (U^n)$ and $||\varphi||_\infty \leq ||T_{\bar{\varphi}}||.$
\end{theorem}

\begin{theorem}
	Let $0<\max (p,q)\leq 1, 1<s<\infty, m,k \in \mathbb{N}, k=\alpha+\frac{1}{p}-\frac{1}{s}+m.$ Then $T_\varphi$ operator is bounded operator from $F^{p,q}_{\alpha,k} (U^n)$ into $H^{s}_{m} (U^n)$ if and only if $\varphi \in H^\infty (U^n)$ and $||\varphi||_\infty \leq ||T_{\bar{\varphi}}||.$
\end{theorem}

The proof of Theorem 1.

Note if $f_r(z_1,...,z_n)=\prod_{i=1}^m$$\frac{1}{1-r_jz_j}, z_j \in U, r_j\in (0,1), j=1,...,n.$ Then

\begin{align*}
	 T_{\bar{h}}(f_r)=\frac{1}{(2\pi i)^n}\int_{T^n}\frac{f_r(\xi_1,...,\xi_n)\overline{h(\xi_1,...,\xi_n)}}{\prod_{j=1}^n(\xi_j-z_j)}d\xi_1...d\xi_n=
     \\=\frac{1}{(2\pi i)^n} \int_{T^n}\frac{\overline{h(\xi_1,...,\xi_n)}d\xi_1...d\xi_n}{\prod_{j=1}^n (\xi_1-z_j)(1-r_j\xi_j)}=\\
	=\frac{(-1)^n}{(2\pi i)^n}\int_{T^n}\frac{h(\xi_1,...,\xi_n)d\bar{\xi}_1...d\bar{\xi}_n}{\prod_{j=1}^n(\bar{\xi}_j-\bar{z}_j)(1-r_j\bar{\xi}_j)}=\\
	=\frac{1}{(2\pi i)^n}\int_{T^n}\frac{h(\xi_1,...,\xi_n)d\xi_1...d\xi_n}{\prod_{j=1}^n(\bar{\xi}_j-\bar{z}_j)\xi^2_j(1-r_j\bar{\xi}_j)}=\\
	=\frac{1}{(2\pi i)^n}\int_{T^n}\frac{h(\xi_1,...,\xi_n)d\xi_1...d\xi_n}{\prod_{j=1}^n(1-\bar{z}_j\xi_j)(\xi_j-r_j)}=\\
	=\frac{\overline{h(r_1,...,r_n)}}{\overline{\prod_{j=1}^n(1-\bar{z}_jr_j)}}=	\frac{\overline{h(r_1,...,r_n)}}{\prod_{j=1}^n(1-z_jr_j)}, z_j \in U^n, j=1,...,m.\\
\end{align*}

Using these equalities we can show one part of our theorem (necessity of condition
on symbol).
The necessity of condition $\varphi \in H^\infty(U^n)$ can be shown in a standard way. Let $T\bar{\varphi}$
be a bounded operator from $F^{p,q}_{\alpha,k} (U^n)$ into $A^{s}_{m,m}(U^n)$ . We consider the following
test function

$$
f_r(z)=\frac{(1-r)^{k+1-\alpha-\frac{1}{p}}}{(1-rz)},
$$

where $(1-rz)=\prod_{j=1}^{n}(1-r_jz_j), r \in I^n, z\in U^n, r_j>1/2, j=1,...,n.$ Note it
is easy to check using well-known Forelly-Rudin type estimates and the chain of
equalities we provided in the start of our proof that the following two estimates are
valid.

\begin{flushleft}
	(A) $||f_r(z)||_{F^{p,q}_{\alpha,k}}\leq C, k>\frac{1}{p}+\alpha-1$ \\
	(B) $||T_{\bar{\varphi}}(f_r)(z)||_{A^{s}_{m,m}}= r|\varphi(r)| ||f_r(z)||_{A^{s}_{m,m}} \geq r|\varphi(r)|C', r \in (C,1).$
\end{flushleft}

Namely we use the following standard equation to get (A), (B)

$$
D^\alpha\left(\frac{1}{1-rz}\right)=\frac{1}{(1-rz)^{\alpha+1}}, \alpha >0, r \in (0,1), z \in U^n,
$$

and the estimates (see [1], [2], [6], [9], [10])

$$
\int_{U^n}\frac{(1-|z|)^tdm_{2n}(z)}{|1-\bar{z}w|^v}\asymp \frac{1}{(1-|w|)^{v-t-2}}, w \in U^n, t>-1, v>t+2,
$$

where $|1-zw|^\alpha=\prod_{j=1}^{n}|1-z_jw_j|^\alpha, z,w \in U^n, \alpha>0,$

$$
\int_{T^n}\frac{dm_{n}(\xi)}{|(1-rz)^\alpha|}\asymp \frac{1}{(1-rR)^{\alpha-1}}, \alpha>1, r,R \in (0,1), z=R\xi, x \in U^n.
$$

Hence from (A) and (B) we obtain directly that $|\varphi(r)| \leq C',  r \in I^n$, from some
positive constant $C'$.

Note for $|\varphi(r)| \leq C'',  z\in U^n$ we must consider $\tilde{f}_r(z)=f_r(e^{-i\theta}z), \theta \in T^n, r \in (0,1), e^{-i\theta}z=(e^{-i\theta_1}z_1, ..., e^{-i\theta_n}z_n)$. The main problem is now to show the
reverse (sufficiency condition symbol).

Let $\varphi \in H^\infty(U^n)$. Let $F(w)=T_{\bar{\varphi}}f(w), R \in (0,1)$. Then using standard duality
arguments it is easy to show that

$$
||F_R(w)||^s_{A^s_{m,m}}=\int_{U^n}|D^mF_R(w)|^s(1-|w|)^{m-1}dm_{2n}(w)=
$$

$$
=\left|\int_{U^n} D^mF_R(w)\overline{G_R(w)}(1-|w|)^{m-1}dm_{2n}(w)\right|,
$$

where $G_R(w) \in L^{s'}_m, ||f||^{s'}_{L^{s'}_m}=\int_{U^n}|f(z)|^{s'}(1-|z|)^{m-1}dm_{2n}(z)< \infty, \frac{1}{s}+\frac{1}{s'}=1,$ and
$||G_R(w)||_{L^{s'}_m}=1$. Note further 
$F_R(w)=()T_{\tilde{\varphi}}f)(Rw)=\frac{1}{(2\pi)^n}\int \frac{\bar{\varphi}(t)f(t)}{1-\bar{t}Rw}dm_n(t)$.

And hence we have

$$
||F_R(w)||^s_{A^s_{m,m}}
=C(n)\left|\int_{U^n}\int_{T^n} \frac{\overline{\varphi(t)}f(t)\overline{G_R(w)}(1-|w|)^{m-1}}{(1-\bar{t}Rw)^{m+1}}dm_{2n}(w)dm_n(t)\right|.
$$

Using Fubinis theorem we note that the following map

$$
S(G_R)=\Psi(\bar{t}R)
=\int_{U^n} \frac{\overline{G_R(w)}(1-|w|)^{m-1}}{(1-\bar{t}Rw)^{m+1}}dm_{2n}(w)=
$$

$$
=\int_{U^n} \frac{\overline{G_R(z)}(1-|z|)^{m-1}}{(1-\bar{t}R\bar{z})^{m+1}}dm_{2n}(\bar{z}), z=\bar{w},
$$

is a Bergman projection map acting from $L^{s'}_m$ into $A^{s'}_m$, $A^{s'}_m=L^{s'}_m\cap H(U^n)$, $\frac{1}{s}+\frac{1}{s'}=1$.  This important known fact have many applications and can be seen
in [1], [6], [8], [9]. Using Littlewood-Paley equality in the unit polydisk with $D^\alpha$
operator namely the following equality

\begin{flushleft}
	(C) $\frac{1}{(2\pi)^n}\int_{T_n}f(rt)g(r\bar{t})dm_n(t)=$ \\
	 $=(\frac{m}{n})^n\prod\limits^{n}_{j=1}(r_j^{-2m})\int_{0}^{r_1}...\int_{0}^{r_n}\int_{T^n}(D^mg(R\xi))(f(R\bar{\xi}))\prod\limits^{n}_{i=1}(r_i^2R_i^2)^{m-1}RdRdm_n(\xi),$
\end{flushleft}
$f,g \in H(U^n), r \in I^n, m \in \mathbb{N} (see [1], [2], [5], [7]-[10]).$

Hence we have using this equality and Holder’s inequality in the unit polydisk

$$
||F_R(w)||_{A^s_{m,m}(U^n)}
=C(n)\left|\int_{T^n}\phi(\bar{t}R)f(t)\overline{\varphi(t)}dm_n(t)\right|, \phi(w) \in A^{s'}_m.
$$

Then we have the following

\begin{align*}
	||F(w)||_{A^s_{m,m}(U^n)}=\lim\limits_{R\to 1}||F_R(w)||_{A^s_{m,m}(U^n)}=\\
	=C(n)\left|\lim\limits_{R\to 1}\int_{T^n}(\phi(\bar{t}R))(f(t))(\overline{\varphi(t)})dm_n(t)\right|=\\
	=C(n)\left|\lim\limits_{R\to 1}\int_{T^n}(\phi(\bar{t}R))f(Rt)(\bar{\varphi}(\bar{t}R))dm_n(t)\right|.
\end{align*}

Using (C) we have

\begin{align*}
	||F(w)||_{A^s_{m,m}(U^n)}\leq \tilde{c}\int_{U^n}|D^kf(\bar{w})|(1-|w|)^{k-1}|\bar{\varphi}(w)||\phi(w)|dm_{2n}(w)\leq\\
	\leq\tilde{\tilde{c}}||\varphi||_\infty\int_{U^n}|D^kf(w)|(1-|w|)^{k-1}|\phi(w)|dm_{2n}(w).
\end{align*}

From this estimate we using Holder’s inequality have now

\begin{align*}
\int_{U^n}|D^kf(w)|(1-|w|)^{k-1}|\phi(w)|dm_{2n}(w)\leq\\
	\leq \left( \int_{U^n}|D^kf(w)|^s(1-|w|)^{(k-1)s-\frac{s}{s'}(m-1)}dm_{2n}(w)\right)^{\frac{1}{s}}\times \\
	\times\left( \int_{U^n}(1-|w|)^{m-1}|\psi(w)|^{s'}dm_{2n}(w)\right)^{\frac{1}{s'}}, \frac{1}{s}+\frac{1}{s'}=1.
\end{align*}

It remains to use the following lemma (see [4]) for $F^{
p,q}_{\alpha,k}$ type analytic function
spaces.

\begin{lemma} Let $f\in F^{
		p,q}_{\alpha,k} (U^n), s>1, k\in \mathbb{N}, 0<p,q\leq s$ and $k=\frac{(\alpha+1/p)s-1}{s}+m(1-1/s).$ Then the following equality is valid
\begin{align*}
	\left( \int_{U^n}|D^kf(w)|^s(1-|w|)^{(k-1)s-(m-1)(1-1/s)}dm_{2n}(w)\right)^{\frac{1}{s}}\leq \\
	\leq c\left(\int_{T^n}\left( \int_{I^n}|D^kf(w)|^q(1-|w|)^{\alpha q-1}dm_{n}|w|\right)^{\frac{p}{q}}dm_n(\xi)\right)^{\frac{1}{p}}.
\end{align*}		
		
\end{lemma}

Theorem is proved.

\textbf{Remark 1.} Note similar proof can be provided for a little bit general situation where
standard $(1-|w|)^\alpha$  
type weights are replaced with
 $w(r), r \in (0,1)$
 type weights
(see for these weights, for example, [4], [9]).

The proof of Theorem 2.

The necessity of the $\varphi \in H^\infty(U^n)$ condition can be shown as in our previous
theorem. It is enough to consider the following test function

$$
f_r(z)=\frac{(1-r)^{k+1-\alpha-1/p}}{1-rz}
$$

$(1-rz)=\prod_{i=1}^{n}(1-r_iz_i), r\in I^n, z \in U^n, k>\alpha+1/p-1$ and note in addition that $||f^r(z)||_{H^s_m}\geq C$ for $k-m=1/p-1/s+\alpha.$ We show now that $\varphi \in H^\infty (U^n) $  condition
is also sufficient.

Let $G(z)=(T\bar{\varphi}f)(z), R \in (0,1).$ Then using standard duality arguments we
have the following equality

\begin{align*}
||D^mG_R||_{H^s(U^n)}=\left(\int_{T^n}|D^mG(Re^{i\varphi})|^sdm_n(\varphi)\right)^{\frac{1}{s}}=\\
=\left|\int_{T^n} D^m G (Re^{i\varphi})\overline{\phi(e^{i\varphi})} dm_n(\varphi)\right|, R\in (0,1), \phi \in L^{s'}(T^n).
\end{align*}

Note then that $G(z)=T_{\bar{\varphi}}(f)(z)=\frac{1}{(2\pi i)^n}\int_{T^n}\frac{f(t)\overline{\varphi(t)}}{\prod_{j=1}^{n}(1-\bar{t}_jz_j)}dm_n(t), \varphi \in H^\infty(U^n)$,   $z=Re^{i\varphi},$  we have that

 $$D^mG(z)=\frac{1}{(2\pi i)^n}\int_{T^n}\frac{f(t)\overline{\varphi(t)}}{(1-\bar{t}z)^{m+1}}dm_n(t), z=Re^{i\varphi}.$$

And hence we have the following

\begin{align*}
	||D^mG||_{H^s(U^n)}=\lim_{R \to 1} ||D^mG_R||_{H^s(U^n)}=\\
	=C(n) \left|\lim_{R \to 1}\int_{T^n} \int_{T^n}\overline{\phi(e^{i\phi})} \frac{f(t)\overline{\varphi(t)}}{(1-\bar{t}Re^{i\psi})^{m+1}}dm_n(t)dm_n(\psi)\right|=\\
	=C(n) \left|\lim_{R \to 1}\int_{T^n}D^m\left( \int_{T^n} \frac{\overline{\varphi(e^{i\psi})}dm_n(\psi)}{(1-\bar{t}Re^{i\psi})}\right)f(t)\overline{\varphi(t)}dm_n(t)\right|.
\end{align*}

Note now that the following operator (see [8], [9])

$$
(S(\phi))(\bar{t}R)=h(\bar{t}R)= \int_{T^n} \frac{\overline{\phi(e^{-i\psi})}dm_n(\psi)}{1-\bar{t}Re^{i\psi}}
$$

is a Riesz projection of $\bar{\phi}$ function and hence $h(\bar{t}R) \in H^{s'}(T^n), \frac{1}{s}+\frac{1}{s'}=1$.
This important fact concerning boundedess of Riesz projection in Hardy spaces in
polydisks can be seen in [8]. Hence by Littlewood-Paley identity (see above (C))
we have the following inequalities

\begin{align*}
	||D^mG||_{H^s(U^n)}=C(n) \left|\lim_{R \to 1}\int_{T^n}f(Rt)\bar{\varphi}(R\bar{t})D^mh(\bar{t}R)dm_n(t)\right|\leq\\
	\leq C_1(n)\int_{U^n}|D^kf(\bar{w})|(1-|w|)^{k-1}|\bar{\varphi}(w)|D^mh(w)dm_{2n}(w)\leq\\
	\leq C_1(n)\int_{U^n}|D^kf(\bar{w})|(1-|w|)^{k-1}|D^mh(\bar{w})|dm_{2n}(w)||\varphi||_\infty,
\end{align*}

$w=\rho\xi, h(w)\in H^{s'}(T^n), \frac{1}{s}+\frac{1}{s'}=1.$

Using Holders inequality we finally have

$$
	||D^mG||_{H^s(U^n)}\leq C_1(n) ||\varphi||_\infty \int_{I^n}M_s(D^kf,\rho)M_{s'}(D^mh,\rho)(1-\rho)^{k-1}\rho d\rho_1...d\rho_n.
$$

Assume now that the following lemma is valid (see the proof below).

\begin{lemma} Let $G\ in H(U^n), 1<s<\infty.$ Then the following estimates are valid
\begin{align*}
	M_s(D^mG,R^2)\leq c(1-R)^{-m}M_s(G,R), R\in I^n, m\in \mathbb{N}\\
	M_s(G,R^2)\leq c_1(1-R)^{1/s-1}M_1(G,R), R\in I^n.\\
\end{align*}
\end{lemma}

\textbf{Remark 2.} Note it is well known that those estimates are valid for $n = 1$ case (see
[2], [5], [6], [9]).
Now using this lemma from last estimates we obtained above, we get the following
estimates for the $||D^mG||_{H^s(U^n)}$.  We have

$$
||D^mG||_{H^s(U^n)}\leq c ||\varphi||_\infty ||h||_{H^{s'}(U^n)} \int_{I^n}M_s(D^kf,\rho)(1-\rho)^{k-m-1}\rho d\rho_1...d\rho_n,
$$

$k>m, k, m \in \mathbb{N}  $. Next, from Lemma B and Lema 2 (see below) and Lemma A we
have the following estimates

\begin{align*}
\int_{U^n}|G(w)|(1-|w|)^{1/p+\alpha-2}dm_{2n}(w)\leq\\
\leq c\left(\int_{T^n}\left(\int_{I^n}(1-R)^{\alpha q-1}|G(R\xi)|^q d R\right)^{\frac{p}{q}}dm_n(\xi)\right)^{\frac{1}{p}}, G \in H(U^n), max(p,q)\leq 1, 0<\alpha<\infty,\\
\int_{I^n}M_s(G,R)(1-R)^\beta RdR\leq  c\int_{I^n}|G(w)|(1-|w|)^{1/s-1+\beta}dm_{2n}(w), G \in H(U^n), \beta>-1, s>1.
\end{align*}

Using this last estimates two estimates we have

$$
||D^mG||_{H^s(U^n)}\leq c ||\varphi||_\infty ||h||_{H^{s'}(U^n)} ||f||_{F^{p,q}_{k,\alpha}}<\infty.
$$

Theorem is proved.

The proof of Lemma A.
The second part follows directly from Lemma 2 (see [11], [12]), the first part is
easy to get using Fubinis theorem and induction.

\begin{lemma} (see [11], [12]) Let f be analytic in $0\leq r_j<|z_j|<R_j, 1\leq j\leq m$ and $f \in C^s$ continuous in closure of this domain. Then for $0\leq p\leq q\leq \infty, \rho \in (r,R)$

$$
||f_\rho||_{H^q}\leq c(m,p,q)\prod_{j=1}^{m}((\rho_j-r_j),(R_j-\rho_j))^{1/q-1/p} \max_{V_j=r_j,R_j\\ j=1,...,m} ||f_v||_{H^p}.$$	
\end{lemma}

\begin{lemma}
(see [9], [10]) Let $0<max(p,q)\leq s < \infty, \alpha>0.$ Then

$$
\left(\int_{U^n}|f(w)|^s(1-|w|)^{s(\alpha+1/p)-2}dm_{2n}(w)\right)^{\frac{1}{s}}\leq c||f||_{F^{p,q}_\alpha}.
$$
\end{lemma}

\textbf{Remark 3.} It is easy to note that the same approach can be used to get criteria
on symbol $\varphi$, for which $(T_{\bar{\varphi}})$ operator is bounded from  $F^{p,q}_{k,\alpha}$ into $X$, where $X$ is a
different from $A^m_s$ and $H^s$ quazinormed subspace of $H(U^n).$

 \end{document}